\documentclass[12pt]{amsart}
\usepackage{amssymb}
\usepackage{amsfonts}
\usepackage{amssymb,latexsym}
\usepackage{enumerate}
\usepackage{mathrsfs}
\usepackage{hyperref}
 \usepackage{color}\makeatletter
\@namedef{subjclassname@2010}{%
  \textup{2010} Mathematics Subject Classification}
\makeatother

  [2002/01/19 v2.2g %
    AMS font definitions%
  ]
\DeclareFontFamily{U}{euf}{}
\DeclareFontShape{U}{euf}{m}{n}{%
  <5><6><7><8><9>gen*eufm%
  <10><10.95><12><14.4><17.28><20.74><24.88>eufm10%
  }{}
\DeclareFontShape{U}{euf}{b}{n}{%
  <5><6><7><8><9>gen*eufb%
  <10><10.95><12><14.4><17.28><20.74><24.88>eufb10%
  }{}

  [2002/01/19 v2.2g %
    AMS font definitions%
  ]
\DeclareFontFamily{U}{msb}{}
\DeclareFontShape{U}{msb}{m}{n}{%
  <5><6><7><8><9>gen*msbm%
  <10><10.95><12><14.4><17.28><20.74><24.88>msbm10%
  }{}

  [2002/01/19 v2.2g %
    AMS font definitions%
  ]
\DeclareFontFamily{U}{msa}{}
\DeclareFontShape{U}{msa}{m}{n}{%
  <5><6><7><8><9>gen*msam%
  <10><10.95><12><14.4><17.28><20.74><24.88>msam10%
  }{}

\newtheorem{theorem}{Theorem}[section]
\newtheorem{lemma}[theorem]{Lemma}
\newtheorem{proposition}[theorem]{Proposition}

\theoremstyle{definition}
\newtheorem{definition}[theorem]{Definition}

\newtheorem{remark}[theorem]{Remark}

\numberwithin{equation}{section} 

\begin{document}

\title[]
{Genus theory and Euclidean ideals for real biquadratic fields}

\begin{abstract}
In this paper, we use the theory of genus fields to study the Euclidean ideals of certain real biquadratic fields $K.$ Comparing with the previous works, our methods yield a new larger family of real biquadratic fields $K$ having Euclidean ideals;
 and the conditions for our family   seem to be more efficient for the computations. Moreover, the previous approaches  mainly focus  on  the case
if $h_K=2$,  while the present approach  can also deal with the general case when $h_K=2^t  (t\geq1)$, where $h_K$ denotes the ideal class number of $K$. In particular, if $h_K\geq 4$, it shows that $H(K)$, the Hilbert class field of $K$, is always non-abelian over $\mathbb{Q}$ for the family of $K$ given in this paper having Euclidean ideals, whereas the previous approaches  always requires that $H(K)$ is abelian over $\mathbb{Q}$ explicitly or implicitly. Finally, some open questions have also been listed for further research.
\end{abstract}

\author{Su Hu}
\address{Department of Mathematics, South China University of Technology, Guangzhou, Guangdong 510640, China}
\address{State Key Laboratory of Cryptology, P. O. Box 5159, Beijing, 100878, China}
\email{mahusu@scut.edu.cn}

\author{Yan Li*}
\address{Department of Applied Mathematics, China Agricultural university, Beijing 100083, China}
\address{State Key Laboratory of Cryptology, P. O. Box 5159, Beijing, 100878, China}
\email{liyan\_00@cau.edu.cn }

\thanks{*Corresponding author}



\maketitle

\vspace*{1ex}

\def\C{\mathbb C_p}
\def\BZ{\mathbb Z}
\def\Z{\mathbb Z_p}
\def\Q{\mathbb Q_p}
\def\nu{{\rm ord}}

\section{Introduction}\label{intro}
An  integral domain $R$ is  Euclidean if it has an Euclidean algorithm, i.e., a function
$$\psi: R \rightarrow \mathbb{N}=\{0,1,2,\ldots\}$$ satisfying the property:
for all $a,b\in R$, if $b\neq 0$, then there exists $q,r\in R$ such that $a=bq+r$ with  $\psi(r) < \psi(b)$.

In 1979 H.W. Lenstra  generalized the concept of Euclidean algorithm
from rings to ideals.
\begin{definition}[{Lenstra \cite{Lenstra}}] Let $R$ be a Dedekind domain and $\mathbb{I}$ be the set of all fractional ideals containing $R$. A fractional ideal $C$ of $R$ is said to be an Euclidean ideal if there exists a function $\psi : \mathbb{I} \rightarrow W$, for some well-ordered set $W$, such that for all $I \in \mathbb{I}$ and all $x \in IC \setminus C$, there exists some $y \in C$ such that
\begin{equation}
\psi((x-y)^{-1}IC) < \psi (I).
\end{equation}
\end{definition}
Clearly, if $C$ is an Euclidean ideal, then every ideal in the ideal class $[C]$ is also Euclidean.
So we can say Euclidean ideal class.

The importance of this new definition is 
that if $R$ has an Euclidean ideal $C$, then the ideal class group of $R$ is cyclic and generated by $C$ (\cite[Theorem 1.6]{Lenstra}) just as the well known fact: Euclidean domain is principal. 


For the case of integer rings $\mathcal{O}_{K}$ of number fields $K$, Lenstra \cite{Lenstra} proved a ``converse" statement holds under the truth of the generalized Riemann hypotheses  (GRH), that is, assuming  GRH,  all the number fields $K$ with cyclic class group,
except the imaginary quadratic fields, have an Euclidean ideal. He also stated that $\mathbb{Q}(\sqrt{-d})$ with $d=1,2,3,5,7,11,15$ are the only quadratic imaginary fields with an Euclidean ideal.





Without assuming GRH, Graves and Ram Murty proved the following general theorem by sieve methods in 2013.

\begin{theorem}[{Graves and Ram Murty~\cite{Graves 2}}]\label{Theorem 4.1} Let $K$ be a number field, with ring of integers $\mathcal{O}_{K}$, group of units  $\mathcal{O}_{K}^*$ and cyclic class group $\mathrm{Cl}(K)$. If its Hilbert class field, $H(K)$, has an abelian Galois group over $\mathbb{Q}$ and if ${\rm rank} (\mathcal{O}^*_{K})\geq 4$, then
$$
\mathrm{Cl}(K)=\langle[C]\rangle\ if\ and\ only\ if\ [C]\ is\ an\ Euclidean\ ideal\ class.
$$
\end{theorem}

They also gave an example: $K=\mathbb{Q}(\sqrt{5}, \sqrt{21}, \sqrt{22})$ with Hilbert class field
$H(K)=\mathbb{Q}(\sqrt{5}, \sqrt{3}, \sqrt{7}, \sqrt{2}, \sqrt{11})$. During the same period, Graves found the first example of totally real number field $K$ with unit rank $3$ having a non-principal Euclidean ideal. 

\begin{theorem}[{Graves \cite{Graves 4}}]The field $\mathbb{Q}(\sqrt{2},\sqrt{35})$ has a non-principal Euclidean ideal.
\end{theorem}


Extending the work of Graves \cite{Graves 4},  Hsu~\cite{Hsu} explicitly constructed a family of real biquadratic fields having a non-principal Euclidean ideal class. Her result reads as follows.
\begin{theorem}[{Hsu \cite{Hsu}}]\label{HsuT} Suppose $K$ is a quartic field of the form $\mathbb{Q}(\sqrt{p_1},\sqrt{q_1q_2})$ with class number $h_K$.  Then $K$ has a non-principal Euclidean ideal class whenever  $h_K = 2$. Here the integers $p_1,q_1,q_2$ are all primes $\geq 29$ and are all congruent to $1$ $\pmod{4}.$
\end{theorem}

In fact, Hsu~\cite{Hsu} also constructed a family of real cyclic quartic fields having such property.
Using PARI,  she computed a list of examples, and conjectured both her classes of number fields are  infinite (See Remark 1.9 of \cite{Hsu}).

Recently, Chattopadhyay and Muthukrishnan~\cite{JNT} further found two other families of biquadratic fields which have a non-principal Euclidean ideal class.
Combing with Hsu's result, they got the following more systematic  result.
\begin{theorem}[{Chattopadhyay and Muthukrishnan \cite[Theorem 1.6]{JNT}}]\label{Theorem 1.7}
	Let $K=\mathbb{Q}(\sqrt{p_1},\sqrt{q_1q_2})$, where $p_1$ is any rational prime and $q_1,q_2 \equiv 1\pmod 4$ are prime numbers. If $h_{K}=2$, the $K$ has an  Euclidean ideal class.
\end{theorem}

All above mentioned results concerning number fields $K$ with Euclidean ideal class relies
on the fact that the Hilbert class field $H(K)$ is abelian over $\mathbb{Q}$, and except
$K=\mathbb{Q}(\sqrt{5}, \sqrt{21}, \sqrt{22})$, all other examples requires the condition $h_K=2$.

In this paper, based on the genus theory of biquadartic fields, especially the works of Yue~\cite{Yue2010}, Bae $\&$  Yue~\cite{Yue2011}, Ouyang $\&$ Zhang~\cite{Ouyang}, we
relax the condition $h_K=2$ to $h_K=2^t$ and explicitly construct a new larger family of biquadratic fields with Euclidean ideal class. 
In particular, if we assume $h_K\geq 4$, then $H(K)$ is always non-abelian over $\mathbb{Q}$. To our knowledge, this is the first known family of real biquadratic fields having Euclidean ideal and non-abelian Hilbert class field in literature. Here are our main results.
\begin{theorem}\label{main1}
	Let $K=\mathbb{Q}(\sqrt{p_{1}},\sqrt{q_{1}q_{2}})$, where  $p_{1},q_{1},q_{2}$ are three distinct primes. Assume $q_1\neq 3$ and $q_2\neq 3$. Let $G(K)$ be the Hilbert genus field of $K$, i.e. the maximal subfield of Hilbert class field H(K) such that $G(K)$ is a multi-quadratic extension of K. Under the conditions:
\begin{enumerate}[i).]	
\item  $h_{\mathbb{Q}(\sqrt{p_{1}})}=1$; $h_{\mathbb{Q}(\sqrt{q_{1}q_{2}})}$ and $h_{\mathbb{Q}(\sqrt{p_{1}{q_{1}{q_{2}}}})}$ are  powers of $2$;
\item $G(K)=\mathbb{Q}(\sqrt{p_{1}},\sqrt{q_{1}},\sqrt{q_{2}})$;
\end{enumerate}
$K$ has an Euclidean ideal class.  Moreover, if $h_K\geq 4$, then $H(K)$ is always non-abelian over $\mathbb{Q}$.
\end{theorem}
\begin{remark}\label{Kurodarem}

By Kuroda's class number formula, the condition $i)$ of Theorem \ref{main1} is equivalent to
$$h_{K}\ is\ a\ power\ of\ 2.$$ Indeed, for any real biquadratic field $K$ with intermediate quadratic fields $k_{0}, k_{1}, k_{2}$, their class numbers satisfy $$h_K=\frac{1}{4}Q(K)h_{0}h_{1}h_{2},\ {\rm where}\
Q(K)=[\mathcal{O}_{K}^*:\mathcal{O}_{k_0}^*\mathcal{O}_{k_1}^*\mathcal{O}_{k_2}^*]$$
is the unit index; $h_{0},h_{1},h_{2}$ are class numbers of $k_{0},k_{1},k_{2}$, respectively.
In general, the index $Q(K)$ can be $1$, $2$, or $4$ (see \cite{Kuroda} or introductions of \cite{Yue2010} and \cite{Yue4rank}).  Moreover, for real quadratic fields $\mathbb{Q}(\sqrt{d})$ with $d$ square-free, it is well-known that the class number of $\mathbb{Q}(\sqrt{d})$  is odd if and only if $d=p,2p_1, p_1p_2$ such that $p, p_1, p_2$ are primes with $p_1, p_2\equiv 3\mod 4$.
	

\end{remark}


Yue~\cite{Yue2010}, Bae and  Yue~\cite{Yue2011}, Ouyang and Zhang~\cite{Ouyang} explicitly described  Hilbert Genus fields $G(K)$ of all real biquadratic fields $K$ that contains a quadratic subfield $k_0$ of odd class number. Based on their works, we explicitly list all real biquadratic fields $K=\mathbb{Q}(\sqrt{p_{1}},\sqrt{q_{1}q_{2}})$ with  $p_{1},q_{1},q_{2}$ three distinct primes, whose Hilbert genus field $G(K)$ is equal to $\mathbb{Q}(\sqrt{p_{1}},\sqrt{q_{1}},\sqrt{q_{2}})$.
This simplifies the condition $ii)$ of Theorem \ref{main1}.

\begin{theorem}\label{genus}
	Let $K=\mathbb{Q}(\sqrt{p_{1}},\sqrt{q_{1}q_{2}})$, where $p_{1}, q_{1}, q_{2}$ are three distinct primes. Let $\varepsilon$ be the fundamental unit of $k=k_0=\mathbb{Q}(\sqrt{p_{1}})$ and $\left(\frac{\bullet}{\bullet}\right)$ be the Legendre symbol. 
Then $G(K)=\mathbb{Q}(\sqrt{p_{1}},\sqrt{q_{1}},\sqrt{q_{2}})$ if and only if $p_{1}, q_{1}, q_{2}$ belongs to one of the following cases.\\
Case 1. $p_1 \equiv 1 \pmod4$, $d=q_{1}q_{2} \equiv 1 \pmod4$
\begin{enumerate}[$\bullet$]
	\item $q_{1}, q_{2} \equiv 1 \pmod4$, $\displaystyle\left(\frac{p_{1}}{q_{1}}\right)=\displaystyle\left(\frac{p_{1}}{q_{2}}\right)=-1$
	\item $q_{1}, q_{2} \equiv 1 \pmod4$, $\displaystyle\left(\frac{p_{1}}{q_{1}}\right)=1, \displaystyle\left(\frac{p_{1}}{q_{2}}\right)=-1, \displaystyle\left(\frac{\varepsilon \mod \mathfrak{q} }{q_{1}}\right)=-1$
where $\mathfrak{q}$ is some (equivalently any) prime ideal of $\mathcal{O}_k$ lying above $q_1$.
	\end{enumerate}


\noindent Case 2. $p_{1} \equiv 1 \pmod4$, $d=2q_{1}$
\begin{enumerate}[$\bullet$]
	\item $p_{1} \equiv 5 \pmod8$, $\displaystyle\left(\frac{p_{1}}{q_{1}}\right)=-1, q_{1} \equiv 1 \pmod4$
	\item $p_{1} \equiv 5 \pmod8$, $\displaystyle\left(\frac{p_{1}}{q_{1}}\right)=1, q_{1} \equiv 1 \pmod4, \displaystyle\left(\frac{\varepsilon \mod \mathfrak{q} }{q_{1}}\right)=-1$ where $\mathfrak{q}$ is some (equivalently any) prime ideal of $\mathcal{O}_k$ lying above $q_1$.
\end{enumerate}


\noindent Case 3. $p_{1} \equiv 1 \pmod4$, $d=q_{1}q_{2} \equiv 3 \pmod4$
\begin{enumerate}[$\bullet$]
	\item $p_{1} \equiv 1 \pmod8$, $\displaystyle\left(\frac{p_1}{q_{1}}\right)=\displaystyle\left(\frac{p_1}{q_{2}}\right)=-1$
	\item $p_{1} \equiv 5 \pmod8$, $\displaystyle\left(\frac{p_1}{q_{1}}\right)=\displaystyle\left(\frac{p_1}{q_{2}}\right)=-1$
	\item $p_{1} \equiv 5 \pmod8$, $\displaystyle\left(\frac{p_1}{q_{1}}\right)=1, \displaystyle\left(\frac{p_1}{q_{2}}\right)=-1, q_{1} \equiv 3 \pmod4 $
\end{enumerate}
Case 4. $p_{1} =2$
\begin{enumerate}[$\bullet$]
	\item $q_{1}, q_{2} \equiv 5 \pmod8$
	\item $q_{1} \equiv 1 \pmod8, q_{2} \equiv 5 \pmod8, q_{1} \not\in  A^{+}:=\left\lbrace a^{2}+32b^{2}|a, b \in\mathbb{Z}\right\rbrace.$
	\item $q_{1} \equiv 5 \pmod8, q_{2} \equiv 3 \pmod8$
	\item $q_{1} \equiv 7 \pmod8, q_{2} \equiv 5 \pmod8$
\end{enumerate}
Case 5. $p_{1} \equiv 3 \pmod4$, $d=q_{1}q_{2} \equiv 1\,\textrm{or}\,3 \pmod4$
\begin{enumerate}[$\bullet$]
	\item $\displaystyle\left(\frac{p_{1}}{q_{1}}\right)=\displaystyle\left(\frac{p_{1}}{q_{2}}\right)=-1$
	\item $\displaystyle\left(\frac{p_{1}}{q_{1}}\right)=1, \displaystyle\left(\frac{p_{1}}{q_{2}}\right)=-1, q_{1} \not\equiv 1 \pmod8 $
	\item $\displaystyle\left(\frac{p_{1}}{q_{1}}\right)=\displaystyle\left(\frac{p_{1}}{q_{2}}\right)=1, q_{1} \not\equiv 1 \pmod8, q_{2} \not\equiv 1 \pmod8, q_{1}q_{2} \equiv 5,7 \pmod8 $
\end{enumerate}
Case 6. $p_{1} \equiv 3 \pmod4$, $d=2q_{1}$
\begin{enumerate}[$\bullet$]
	\item $\displaystyle\left(\frac{p_{1}}{q_{1}}\right)=-1$
	\item $\displaystyle\left(\frac{p_{1}}{q_{1}}\right)=1, q_{1} \not\equiv 1 \pmod8$.
\end{enumerate}
\end{theorem}
\begin{remark}In all cases of Theorem \ref{genus}, the conditions are very easy to check. Most times, one only needs to compute the Legendre symbols. In the second subcase of of \emph{Case 1} and \emph{Case 2}, one needs to compute the fundamental unit of $\mathbb{Q}(\sqrt{p_1})$ additionally, while in the second subcase of \emph{Case 4}, one needs to check that $q_1$ can not be represented by the binary integral quadratic form $a^2+32b^2$.
  So combining Theorem \ref{main1} and Theorem \ref{genus}, it is easy to produce examples of real biquadratic fields having Euclidean ideals by computer software, since we mainly need to compute the class number $h_{k_0}, h_{k_1},  h_{k_2}$ of quadratic subfields of $K$ and Legendre symbols. Note that previous methods needs to compute $h_K$, which seems to cost more labour.
\end{remark}
In the special case: $q_1\equiv q_2\equiv1\pmod{4}$, we get a complete answer under the condition $h(K)=2^t$, which also extends previous results of Hsu, {Chattopadhyay and Muthukrishnan \cite[Theorem 1.6]{JNT}} from $h(K)=2$ to $h(K)=2^t$.

\begin{theorem}\label{main3}
	Let $K=\mathbb{Q}(\sqrt{p_{1}},\sqrt{q_{1}q_{2}})$, where  $p_{1},q_{1},q_{2}$ are three distinct primes such that $q_1\equiv q_2\equiv1\pmod{4}$.  Assume that $h_K$ is a power of $2$, or equivalently
$$h_{\mathbb{Q}(\sqrt{p_{1}})}=1; \ h_{\mathbb{Q}(\sqrt{q_{1}q_{2}})}\ \mathrm{ and}\  h_{\mathbb{Q}(\sqrt{p_{1}{q_{1}{q_{2}}}})}\ \mathrm{are\   powers\  of}\ 2.$$
  Then $K$ has an  Euclidean ideal class \textbf{if and only if}
$$G(K)=\mathbb{Q}(\sqrt{p_{1}},\sqrt{q_{1}},\sqrt{q_{2}}),$$
where $G(K)$ is the Hilbert genus field of $K$; or equivalently, 
$ q_{1}\equiv q_{2}\equiv 1\pmod{4}$ and
 $p_{1}, q_{1}, q_{2}$ belongs to one of the following cases:

Case 1. $p_1 \equiv 1 \pmod4$
\begin{enumerate}[$\bullet$]
	\item $\displaystyle\left(\frac{p_{1}}{q_{1}}\right)=\displaystyle\left(\frac{p_{1}}{q_{2}}\right)=-1$
	\item  $\displaystyle\left(\frac{p_{1}}{q_{1}}\right)=1, \displaystyle\left(\frac{p_{1}}{q_{2}}\right)=-1, \displaystyle\left(\frac{\varepsilon \mod \mathfrak{q} }{q_{1}}\right)=-1$
where $\varepsilon$ is the fundamental unit of $\mathbb{Q}(\sqrt{p_{1}})$  and $\mathfrak{q}$ is some (equivalently any) prime ideal of $\mathcal{O}_k$ lying above $q_1$.
	\end{enumerate}

Case 2. $p_{1} =2$
\begin{enumerate}[$\bullet$]
	\item $q_{1}, q_{2} \equiv 5 \pmod8$
	\item $q_{1} \equiv 1 \pmod8, q_{2} \equiv 5 \pmod8, q_{1} \not\in  A^{+}:=\left\lbrace a^{2}+32b^{2}|a, b \in\mathbb{Z}\right\rbrace.$
\end{enumerate}

Case 3. $p_{1} \equiv 3 \pmod4$
\begin{enumerate}[$\bullet$]
	\item $\displaystyle\left(\frac{p_{1}}{q_{1}}\right)=\displaystyle\left(\frac{p_{1}}{q_{2}}\right)=-1$
	\item $\displaystyle\left(\frac{p_{1}}{q_{1}}\right)=1, \displaystyle\left(\frac{p_{1}}{q_{2}}\right)=-1, q_{1} \equiv 5 \pmod8 $.
\end{enumerate}
\end{theorem}

Using PARI, Hsu calculated class number $h_K$ for many real biquadratic fields $K=\mathbb{Q}(\sqrt{p_1}, \sqrt{q_1q_2})$ with distinct primes $p_1\equiv q_1\equiv q_2\equiv1\pmod{4}$ and listed them in Table 1 of \cite{Hsu}. If $h_K=2$, then $K$ has an Euclidean ideal by her Theorem(see Theorem \ref{HsuT}). Thus she found many examples of $K$ with Euclidean ideal. But if $h_K>2$, she can not decide it. By our Theorem \ref{main3}, 
we can fill the blanks in her Table whenever $h_K> 2$ is a power of $2$. We list them in Table \ref{tab1}.

By Theorem \ref{main1}, these eighteen fields $K$ with ``Y" in Table \ref{tab1} have non-abelian Hilbert class fields $H(K)$ 
and Euclidean ideals, which, to our knowledge,  is the first time to appear in the literature.

The rest of paper is organized as follows. In section 2, we give some necessary preliminary for the proofs, which includes some basic facts from class field theory, analytic tools of Graves' growth results from \cite{Graves 3} and \cite{Graves 4}, and a brief review of the works of Yue~\cite{Yue2010}, Bae $\&$  Yue~\cite{Yue2011}, Ouyang $\&$ Zhang~\cite{Ouyang} on Hilbert genus fields of biquadratic fields. In section 3, we prove our main results. Finally, in section 4, we list some related open questions on Euclidean ideals.


\begin{table}[htbp]\label{tab1}
  \centering
  \caption{Examples of biquadratic fields $K$ having Euclidean ideals and non-abelian Hilbert Class field}
    \begin{tabular}{ccccccc}
    $(p_1,q_1,q_2)$ & $h_K$     & $\displaystyle\left(\frac{p_1}{q_1}\right)$     & $\displaystyle\left(\frac{p_1}{q_2}\right)$    & $\varepsilon$    & $\displaystyle\left(\frac{\varepsilon \mod \mathfrak{q} }{q_{1}}\right)$       & Euclidean \\
    (29,53,37) & 16    & 1     & $-1$    & $(5+\sqrt{29})/2$ & $-1$      & Y \\
    (29,37,97) & 4     & $-1$    & $-1$    & $(5+\sqrt{29})/2$ &            & Y \\
    (29,41,61) & 4     & $-1$    & $-1$    & $(5+\sqrt{29})/2$ &              & Y \\
    (29,41,89) & 4     & $-1$    & $-1$    & $(5+\sqrt{29})/2$ &             & Y \\
    (29,53,89) & 4     & 1     & $-1$    & $(5+\sqrt{29})/2$ & $-1$      & Y \\
    (29,53,97) & 4     & 1     & $-1$    & $(5+\sqrt{29})/2$ & $-1$     & Y \\
    (37,53,29) & 16    & 1     & $-1$    & $6+\sqrt{37}$    & $-1$      & Y \\
    (37,29,97) & 4     & $-1$    & $-1$    & $6+\sqrt{37}$    &             & Y \\
    (37,41,53) & 4     & 1     & 1     & $6+\sqrt{37}$    &       & N \\
    (37,41,61) & 4     & 1     & $-1$    & $6+\sqrt{37}$    & $-1$      & Y \\
    (37,53,73) & 4     & 1     & 1     & $6+\sqrt{37}$    &       & N \\
    (37,53,89) & 4     & 1     & $-1$    & $6+\sqrt{37}$    & $-1$      & Y \\
    (37,53,97) & 4     & 1     & $-1$    & $6+\sqrt{37}$    & $-1$      & Y \\
    (37,73,61) & 4     & 1     & $-1$    & $6+\sqrt{37}$    & $-1$      & Y \\
    (37,73,89) & 8     & 1     & $-1$    & $6+\sqrt{37}$    & $-1$      & Y \\
    (37,73,97) & 4     & 1     & $-1$    & $6+\sqrt{37}$    & $-1$      & Y \\
    (41,61,29) & 4     & 1     & $-1$    & $32+5\sqrt{41}$  & 1  & N \\
    (41,29,89) & 4     & $-1$    & $-1$    & $32+5\sqrt{41}$  &             & Y \\
    (41,37,53) & 4     & 1     & $-1$    & $32+5\sqrt{41}$  & $-1$  & Y \\
    (41,37,61) & 4     & 1     & 1     & $32+5\sqrt{41}$  &       & N \\
    (41,61,53) & 4     & 1     & $-1$    & $32+5\sqrt{41}$  & 1  & N \\
    (41,61,73) & 8     & 1     & 1     & $32+5\sqrt{41}$  &       & N \\
    (41,61,89) & 4     & 1     & $-1$    & $32+5\sqrt{41}$  & 1  & N \\
    (41,61,97) & 4     & 1     & $-1$    & $32+5\sqrt{41}$  & 1  & N \\
    (41,73,89) & 8     & 1     & $-1$    & $32+5\sqrt{41}$  & $-1$  & Y \\
    (41,73,97) & 4     & 1     & $-1$    & $32+5\sqrt{41}$  & $-1$  & Y \\
    \end{tabular}%
  \label{tab:addlabel}%
\end{table}%


\section{Preliminaries}
\subsection{Artin symbols for abelian extensions}
In this subsection, we review some basic facts of Artin symbols for abelian extensions. The basic reference is Lang \cite[p. 197-200]{LangA}.

Let $K/k$ be an abelian extension of number fields and $\mathfrak{p}$ be a prime of $k$ unramified in $K$. There is an unique element $\sigma\in \textrm{Gal}(K/k)$ satisfying
$$\sigma(\alpha)\equiv \alpha^{N\mathfrak{p}}\pmod{\mathfrak{P}},\ \forall \alpha\in \mathcal{O}_K
$$
where $\mathfrak{P}$ is any prime of $K$ lying above $\mathfrak{p}$ and $N \mathfrak{p}=\#\mathcal{O}_k/\mathfrak{p}$ is the norm of $\mathfrak{p}$. This element $\sigma$ depends only on $\mathfrak{p}$, is denoted by $(\mathfrak{p},K/k)$, and will be called the Artin symbol of $\mathfrak{p}$.

The Artin symbol can be extended multiplicatively to $\mathfrak{a}$, fractional ideals of $k$ unramified in $K$, and still denoted by $(\mathfrak{a}, K/k)$. This is also called the reciprocity law map or the Artin map.

The Artin symbol satisfies the following properties:
\begin{description}
\item[\textbf{A1}] Let $\sigma: K\rightarrow \sigma K$ be an isomorphism (not necessarily to the identity on $k$). Then
                         $$(\sigma\mathfrak{a}, \sigma K/\sigma k)=\sigma (\mathfrak{a}, K/k)\sigma^{-1}. $$
    \item[\textbf{A2}] Let $K'\supset K \supset k$ be a bigger abelian extension. Then
    $$\mathrm{res}_K(\mathfrak{a}, K'/k)=(\mathfrak{a}, K/k)$$
    where the prime ideal factors of $\mathfrak{a}$ should be unramified in $K'$.
    \item[\textbf{A3}] Let $K\supset k' \supset k$ be an intermediate field. Then
$$(\mathfrak{b}, K/k')=(N_{k'/k}\mathfrak{b}, K/k)$$
where for each prime ideal $\mathfrak{q}$ dividing $\mathfrak{b}$ in $k'$, $\mathfrak{p}=\mathcal{O}_k\cap \mathfrak{q}$ should be unramified in $K$. In particular
$$(\mathfrak{q}, K/k')=(\mathfrak{p}, K/k)^f$$
where $f=[\mathcal{O}_{k'}/\mathfrak{q}:\mathcal{O}_{k}/\mathfrak{p}]$ is the residue class degree.
\end{description}

The celebrated Chebotarev density theorem (see Lang \cite[p.169]{LangA}) implies that for each element $\sigma\in \textrm{Gal}(K/k)$, there are \textbf{infinitely many} prime ideals $\mathfrak{p}$ such that
$$(\mathfrak{p}, K/k)=\sigma\ \mathrm{and}\  N\mathfrak{p}=p\ \mathrm{is\ a\ prime}.
$$Therefore the Artin map is surjective.

In the quadratic field case, the Artin symbol equals to the Lengendre symbol, i.e.,
$$(p, \mathbb{Q}(\sqrt{d})/\mathbb{Q})\sqrt{d}=\left(\frac{d}{p}\right)\sqrt{d}
$$
where  $d$ is square-free and $p\nmid d$ is an odd prime.
\subsection{Hilbert symbol}
In this subsection, we review some basic properties of Hilbert symbol for our need. The material here is covered in \cite[Chapter V, \S3 and Chapter VI, \S 8]{Neukirch}.

Let $K$ be a number field, $\mathfrak{p}$ a prime ideal of $K$ with $\mathfrak{p}\cap \mathbb{Z}=p\mathbb{Z}$. Let $K_{\mathfrak{p}}$ be the completion of $K$ at $\mathfrak{p}$. For $a, b\in K_{\mathfrak{p}}^*$, the Hilbert symbol $\left(\frac{a,b}{\mathfrak{p}}\right)\in \mu_2=\{\pm1\}$, is defined by
\begin{equation}\label{Hilbertsymbol}
\left(\frac{a,b}{\mathfrak{p}}\right)=1\Leftrightarrow ax^2+by^2-z^2=0\ \mathrm{has\ a\ nontrivial\ solution\ in}\ K_{\mathfrak{p}}. \end{equation}
The Hilbert symbol defines a non-degenerate bilinear pairing:
$$\left(\frac{\ ,\ }{\mathfrak{p}}\right):\ K_{\mathfrak{p}}^*/K_{\mathfrak{p}}^{*2} \times K_{\mathfrak{p}}^*/K_{\mathfrak{p}}^{*2}\ \rightarrow\ \mu_2
$$
which has the following fundamental properties:
\begin{enumerate}
  \item $\left(\frac{aa',b }{\mathfrak{p}}\right)=\left(\frac{a,b }{\mathfrak{p}}\right)\left(\frac{a',b }{\mathfrak{p}}\right)$
  \item $\left(\frac{a,b }{\mathfrak{p}}\right)=\left(\frac{b,a }{\mathfrak{p}}\right)$
  \item $\left(\frac{a,-a }{\mathfrak{p}}\right)=\left(\frac{a,1-a }{\mathfrak{p}}\right)=1$
  \item $\left(\frac{a,b }{\mathfrak{p}}\right)=1\Leftrightarrow$ $a$ is a norm from the extension $K_{\mathfrak{p}}(\sqrt{b})/K_{\mathfrak{p}}$.
  \item If $\left(\frac{a,b }{\mathfrak{p}}\right)=1$ for all $b\in K_{\mathfrak{p}}^*$, then $a\in K_{\mathfrak{p}}^{*2}$.
\end{enumerate}

Let $\pi$ be a prime element of $K_{\mathfrak{ p}}$, i.e., $\pi$ is the generator of the maximal ideal of integer ring $\mathcal{O}_{\mathfrak{ p}}$. Let $q=\# \mathcal{O}_{\mathfrak{ p}}/\pi \mathcal{O}_{\mathfrak{ p}}$. The following proposition is useful for computing Hilbert symbol.
\begin{proposition}\label{odd}
  Assume $2\nmid p$. For $u, v\in \mathcal{O}_{\mathfrak{ p}}^*$,
  $$\left(\frac{u,v }{\mathfrak{p}}\right)=1\ \mathrm{and}\ \left(\frac{\pi,u }{\mathfrak{p}}\right)=\left(\frac{u}{\mathfrak{p}}\right)$$
  where $\left(\frac{u}{\mathfrak{p}}\right)\equiv u^{\frac{q-1}{2}}\pmod{\pi \mathcal{O}_{\mathfrak{p}}}$ is the quadratic residue symbol modulo $\mathfrak{p}\mathcal{O}_{\mathfrak{ p}}$.
\end{proposition}
For the case of $K_{\mathfrak{p}}=\mathbb{Q}_p$ the $p$-adic field, the Hilbert symbol can be explicitly determined by the following way.

If $p\neq 2$, then
\begin{equation}\label{Qp}
  \left(\frac{p,p}{p}\right)=(-1)^{\frac{p-1}{2}},\ \left(\frac{p,u}{p}\right)=\left(\frac{u}{p}\right),\ \left(\frac{u,v}{p}\right)=1;
\end{equation}
otherwise if $p=2$, then
\begin{equation}\label{Q2}
\left(\frac{2,2}{2}\right)=1,\  \left(\frac{2,u}{2}\right)=(-1)^{\frac{u^2-1}{8}},\ \left(\frac{u,v}{2}\right)=(-1)^{\frac{(u-1)}{2}\frac{(v-1)}{2}}
\end{equation}
where $u, v\in \mathbb{Z}_{p}^{*}$ are units.

For the archimedean place $\mathfrak{p}$, i.e., $\mathfrak{p}$ is a real embedding of $K$, in this case, $K_{\mathfrak{p}}=\mathbb{R}$, or a pair of conjugate complex embedding of $K$, in this case $K_{\mathfrak{p}}=\mathbb{C}$, one can still define the Hilbert symbol $\left(\frac{a,b}{\mathfrak{p}}\right)$ for $a, b\in K_{\mathfrak{p}}^{*}$ by \eqref{Hilbertsymbol}.

These symbols all fit together in the following product formula.
\begin{proposition}\label{product}
  For $a, b\in K^*$, one has
  $$\prod_{\mathfrak{p}}\left(\frac{a,b}{\mathfrak{p}}\right)=1.
  $$
\end{proposition}
\begin{remark}
  In the later use, we will compute the Hilbert symbol for $K$ being a quadratic field. For $p\neq 2$, this can be handled by Proposition \ref{odd}. For $p=2$, if $p$ splits in $K$, then it can be handled by \eqref{Q2}. Otherwise, there is only one prime $\mathfrak{p}$ lying above $2$. This can be handled by the product formula.
\end{remark}

\subsection{Conductors of abelian fields}	
Let $K$ be an abelian field. The celebrated Kronecker-Weber theorem states that $K\subset \mathbb{Q}(\zeta_m)$, $m$th cyclotomic field. The smallest such integer $m$ is called the conductor of $K$, which is denoted by $f(K)$. The following two lemmas on conductors are sufficed for our purpose.

\begin{lemma}[{\cite[Proposition 1]{JNT}}]\label{Proposition1}The conductor $f(K)$ of the real quadratic field $K=\mathbb{Q}(\sqrt{m})$ is
$$
f(K)=\left\{\begin{array}{ll}
m ~~&if~~m\equiv 1(\textrm{mod}~4)\\
4m~~&if~~m\equiv 2,3(\textrm{mod}~4)
\end{array}\right.
$$
where $m>0$ is a squarefree integer.
\end{lemma}


\begin{lemma}[{\cite[Lemma 3]{JNT}}]\label{Lemma 3}
Let $L$ be an abelian number field and let $K_1$, $K_2$ be two subfields of $L$ such that $L$ is the compositum of $K_1$ and $K_2$. Then $f(L)=lcm(f(K_1), f(K_2))$, where $f(K_1), f(K_2)$ and $f(L)$ are the conductors of $K_1$, $K_2$ and $L$, respectively.
\end{lemma}

Lemma \ref{Proposition1} is in fact a special case of Conductor-Discriminant Formula for general abelian fields (See Theorem 3.11 of Washington \cite[p.28]{Washington}).

\subsection{Growth results of Graves}

The following two theorems due to Graves are commonly used analytic tools for proving the existence of Euclidean ideals.
\begin{theorem}[{Graves \cite{Graves 3}}]\label{Theorem 4.1} Suppose that $K$ is a number field such that $|\mathcal{O}_{K}^*|=\infty$, and that $C$ is a non-zero ideal of $\mathcal{O}_K$. If $[C]$ generates the class group of $K$ and
\begin{equation*}
\left| \left\{ \textrm{prime ideal } \wp \subset \mathcal{O}_{K} : N(\wp)\leq x, [\wp]=[C], \pi_{\wp} \mbox{ is onto } \right\} \right|   \gg  \frac{x}{(\log x)^2},
\end{equation*}
where $\pi_{\wp}$ is the canonical map from $\mathcal{O}_K^*$ to $(\mathcal{O}_K/\wp)^*$, then $[C]$ is an Euclidean ideal class.
\end{theorem}
The following result is stated in \cite{Graves 4} directly by utilizing analytical results from \cite{Heath} and \cite{Narkiewicz}.

\begin{theorem}[{Graves~\cite{Graves 4}}]\label{Theorem 1.2}Let $K$ be a totally real number field with conductor $f(K)$ and let $\{e_1, e_2, e_3\}$ be a multiplicatively independent set contained in $\mathcal{O}_{K} ^*$. If $l=lcm(16,f(K))$, and if $gcd(u,l)=gcd(\frac{u-1}{2},l)=1$ for some integer $u$, then
	\begin{align*}
	\left|\left\{ \substack{\textrm{primes }~ \wp \subseteq \mathcal{O}_K \\ \textrm{of degree one}} : N({\wp}) \equiv u \ (mod \ l), N(\wp)\leq x, \langle -1,e_{i}\rangle \twoheadrightarrow (\mathcal{O}_{K}/{\wp})^* \right\}\right| \gg \frac{x}{(\log x)^2},
	\end{align*}
	for at least one $i$.
\end{theorem}
\subsection{Genus fields}
The genus theory has its root in Gauss' classical work on class group of binary quadratic forms. Roughly speaking, genus theory treats the relatively easy part of ideal class groups and class fields. It has been studied and developed by many authors since 1950's. For a more recent development,  see \cite[Chapter IV, \S4]{Gras} for example.

For a number field $K$, the Hilbert genus field of $K$ is defined to be the subfield $G(K)$ of the Hilbert class field $H(K)$ invariant under Gal$(H(K)/K)^{2}$, i.e. the maximal subfield of Hilbert class field $H(K)$ such that $G(K)$ is a multi-quadratic extension of $K$.
From the Artin's reciprocity map, \begin{equation}
\textrm{Cl}(K)\cong \textrm{Gal}(H(K)/K)
\end{equation}
and
\begin{equation}\label{2.2}
\textrm{Cl}(K)/\textrm{Cl}(K)^{2}\cong \textrm{Gal}(G(K)/K).
\end{equation}

Let $K=\mathbb{Q}(\sqrt{p},\sqrt{d})$ be a real quadratic field with $p\equiv 1 (\textrm{mod}~4)$ or $p=2$, and $d$ is
a squarefree positive integer. In 2010, Yue~\cite{Yue2010} constructed the genus field $G(K)$ of $K$ explicitly if $p\equiv 1 (\textrm{mod}~4)$
and $d\equiv 3 (\textrm{mod}~4).$ Then, in 2011, Bae and Yue~\cite{Yue2011} completed Yue's work in the remaining cases.

In 2015, Ouyang and Zhang~\cite{Ouyang} described the genus field $\textrm{G}(K)$ of $K=\mathbb{Q}(\sqrt{\delta},\sqrt{d})$ explicitly, if $\delta=p, 2p$ or $p_{1}p_{2}$
where $p, p_{1}$ and $p_{2}$ are primes congruent to 3 modulo 4, and $d$ is any squarefree positive integer. The works of Yue, Bae and Yue, Ouyang and Zhang
completed the construction of the genus field of real biquadartic fields $K=k_{0}(\sqrt{d})$ such that $k_{0}=\mathbb{Q}(\sqrt{\delta})$ has an odd class number.

For convenience of readers, we give a short review of their works.

For a finite abelian group $A$, let $r_2(A)$ be the $2$-rank of $A$, i.e., the dimension of $A/2A$ or $A[2]$ over  $\mathbb{Z}/2\mathbb{Z}$.

Let $K$ be a real biquadratic field containing a quadratic field $k=k_0$ with odd class number (In our case, $k=\mathbb{Q}(\sqrt{p_1})$ and $K=\mathbb{Q}(\sqrt{p_1},\sqrt{q_1q_2})$). Then
\begin{equation}\label{2dim}
r_2(\textrm{Cl}(K))=s-1-r_2(\mathcal{O}^*_{k}/\mathcal{O}^*_{k}\cap N_{K/k}K^*)
\end{equation}
where $s$ is the number of finite primes of $k$ ramified in $K$ (For a proof, see Proposition 1.1 of \cite{Ouyang}). The number $s$ is relatively easy to compute. The hard part is compute $r_2(\mathcal{O}^*_{k}/\mathcal{O}^*_{k}\cap N_{K/k}K^*)$, which is theoretically based on
\begin{enumerate}
  \item (Hasse Norm Theorem) An element $x\in k^*$ is a norm of $K^*$ if and only if it is a norm locally everywhere, i.e., a norm in every completion $K_{\mathfrak{P}}/k_{\mathfrak{p}}$ ($\mathfrak{P}|\mathfrak{p}$).
  \item (Hilbert symbol) An element $x\in k_{\mathfrak{p}}^*$ is a norm of $K_{\mathfrak{P}}=k_{\mathfrak{p}}(\sqrt{d})$ $\Leftrightarrow$ the Hilbert symbol $\displaystyle\left(\frac{x,d}{\mathfrak{p}}\right)=1$.
\end{enumerate}
Note that in the case $K=\mathbb{Q}(\sqrt{p_1},\sqrt{d})$ with $d>0$ squarefree, $r_2(\mathcal{O}^*_{k}/\mathcal{O}^*_{k}\cap N_{K/k}K^*)$ has an explicit expression except the case that $p_1\equiv 1\pmod{4}$ and $d\equiv 1,2 \pmod{4}$, which was treated in \cite{Yue2011}.

From $r_2(\textrm{Cl}(K))$, one knows the degree of $G(K)/K$. Then they explicitly constructed $K^{*2}\subset \Delta \subset K^{*}$ such that $r_2(\Delta/K^{*2})=r_2(\textrm{Cl}(K))$ and $K(\sqrt{\Delta})/K$ is an unramified extension of $K$. Hence $G(K)=K(\sqrt{\Delta})$. The unramified property of $K(\sqrt{\Delta})/K$ can be checked by passing to local fields. The related extensions of $\mathbb{Q}_2$ is summarized in subsection $2.1$ of \cite{Ouyang}, which can be used to check unramifed property of primes above $2$.

Their computational results show that $G(K)$ is always a multi-quadratic extension of $k$. We do not know how to prove this directly. The $G(K)$ is not necessarily a multi-quadratic extension of $\mathbb{Q}$. For instance:
$$G(K)=\mathbb{Q}(\sqrt{5}, \sqrt{3}, \sqrt{29}, \sqrt{7+2\sqrt{5}})\ \mathrm{ for}\ K=\mathbb{Q}(\sqrt{5},\sqrt{3\cdot29})
$$
(see Example 2.2 of \cite{Yue2010}).
\section{Proof of main results}
\noindent\textbf{Proof of Theorem \ref{main1}:}
By \eqref{2.2} and condition $ii)$,
	$$
	\textrm{Cl}(K)/\textrm{Cl}(K)^{2}\cong  \textrm{Gal}(\mathbb{Q}(\sqrt{p_{1}},\sqrt{q_{1}},\sqrt{q_{2}})/\mathbb{Q}(\sqrt{q_{1}},\sqrt{q_{1}q_{2}}))\cong \mathbb{Z}/2\mathbb{Z}
	$$
Combining condition $i)$, $\textrm{Cl}(K)$ is a cyclic group of order $2^{t}(t\geqslant 1)$.

The proof of $\textrm{Cl}(K)$ having an Euclidean ideal, is based on Theorem \ref{Theorem 4.1} and Theorem \ref{Theorem 1.2}.

Firstly, to apply Theorem \ref{Theorem 4.1},  we characterize the following set:

\begin{equation}\label{3.1}
S=\left\lbrace \textrm{primes}~ \,p:\  2\nmid p, \,p\mathcal{O}_{K}=\prod_{i=1}^{4}\wp_{i}, \,[\wp_{i}] \,\textrm{generates}\, \textrm{Cl}(K)\right\rbrace .\end{equation}
For such prime $p$, $p$ splits completely in $\mathcal{O}_{K}$. Therefore
\begin{equation}\label{3.2}
\displaystyle\left(\frac{p_{1}}{p}\right)=1,\displaystyle\left(\frac {q_{1}q_{2}}{p}\right)=1
\end{equation}
Denote one of $\wp_{i}$ by $\wp$, since $\textrm{Cl}(K)$ is cyclic of order $2^{t}$. $[\wp]$ generates $\textrm{Cl}(K)$ if and only if the image of $[\wp]$ in $\textrm{Cl}(K)/\textrm{Cl}(K)^{2}$ is nontrivial. This can be checked by Artin reciprocity map as follows. Considering
$$\textrm{Cl}(K)/\textrm{Cl}(K)^{2}\rightarrow \textrm{Gal}(G(K)/K)\hookrightarrow \textrm{Gal}(G(K)/\mathbb{Q})$$
$$[\wp] \,\textrm{mod} \,\textrm{Cl}(K)^{2}\mapsto (\wp,G(K)/K)\mapsto (p,G(K)/\mathbb{Q}),$$
as $G(K)=\mathbb{Q}(\sqrt{p_{1}},\sqrt{q_{1}},\sqrt{q_{2}})$, we have
$$\textrm{Gal}(G(K)/\mathbb{Q})\cong  \mu_{2}\times \mu_{2}\times \mu_{2}$$
$$(p,G(K)/\mathbb{Q})\mapsto \left( \displaystyle\left(\frac{p_{1}}{p}\right),\displaystyle\left(\frac{q_{1}}{p}\right),\displaystyle\left(\frac{q_{2}}{p}\right)\right). $$
Combining \eqref{3.1}, $[\wp]$ generates $\textrm{Cl}(K)$ if and only if

 $$\displaystyle\left(\frac{p_{1}}{p}\right)=1, \displaystyle\left(\frac{q_{1}}{p}\right)=-1, \displaystyle\left(\frac{q_{2}}{p}\right)=-1.$$

Therefore
\begin{equation}\label{3.3}
S=\left\lbrace \textrm{primes}~ \,p: \,2\nmid p,\ \displaystyle\left(\frac{p_{1}}{p}\right)=1, \displaystyle\left(\frac{q_{1}}{p}\right)=-1,\displaystyle\left(\frac{q_{2}}{p}\right)=-1 \right\rbrace. \end{equation}

Secondly, we apply Theorem \ref{Theorem 1.2} to show there are enough primes $p$ in $S$, such that $\pi_{\wp}$ is onto. This needs to check $S$ contains an arithmetic progression $p \equiv u \pmod l$, for some $u$ satisfying
\begin{equation}\label{3.31} \textrm{gcd}(u,l)=\textrm{gcd}\left(\frac{u-1}{2},l\right)=1,
\end{equation} where $l =\textrm{lcm}(16,f(K))$.

 By Lemma~\ref{Proposition1} and Lemma~\ref{Lemma 3}, we have $16|l, p_{1}q_{1}q_{2}|l$ and the prime factors of $l$ are $2, p_{1}, q_{1}, q_{2}$.

 If $p_{1}, q_{1}, q_{2}$ are odd, then by the quadratic reciprocity law, the odd primes $p\in S$ is characterized by

\begin{equation}\label{3.4}
\begin{split}\displaystyle\left(\frac{p}{p_{1}}\right)&=(-1)^{\frac{(p_{1}-1)(p-1)}{4}}\\
\displaystyle\left(\frac{p}{q_{1}}\right)&=(-1)^{\frac{(q_{1}-1)(p-1)}{4}+1}\\
\displaystyle\left(\frac{p}{q_{2}}\right)&=(-1)^{\frac{(q_{2}-1)(p-1)}{4}+1}
\end{split}
\end{equation}

Assume  that  $q_1\neq3$ and $q_2\neq3$. If $p_{1}\neq3$, there are at least two quadratic residues modulo each prime $p_{1}, q_{1}, q_{2}$. Therefore, by the Chinese remainder theorem,  we can choose some integer $u$ satisfying \eqref{3.4} (with $p$ replaced by $u$ in \eqref{3.4}) and
\begin{equation}\label{3.41}\begin{aligned} &u \equiv -1 \pmod {4},\ \  u \not\equiv 1 \pmod {p_{1}}, \\&u \not\equiv 1\  \pmod {q_{1}},\ \   u \not\equiv 1 \pmod {q_{2}}. \end{aligned}\end{equation} Otherwise if $p_{1}=3$, we can still choose $u\equiv -1\pmod4$ satisfying equation \eqref{3.4}, which implies that $\displaystyle\left(\frac{u}{3}\right)=-1$, and hence $u \not\equiv 1\pmod3$. So in this case, we can still let $u$  satisfy (\ref{3.41}) and \eqref{3.4}. Since the prime factors of $l=\textrm{lcm}(16,f(K))$ are $2, p_{1}, q_{1}, q_{2}$, the condition \eqref{3.31} automatically holds for such choice of $u$.


Thus if $p \equiv u\pmod l$, we have
$$2\nmid p,\ p  \equiv u \pmod {p_{1}},\  p  \equiv u \pmod {q_{1}},\  p  \equiv u \pmod {q_{2}}.$$
By the choice of $u$, $p$ satisfies \eqref{3.4}. Since $S$ is characterized by (\ref{3.4}), we have $p\in S$.

The case that $p_{1}, q_{1}, q_{2}$ contains 2, can be handled similarly by using the fact
$$\displaystyle\left(\frac{2}{p}\right)= 1 \Longleftrightarrow p \equiv\pm 1\pmod 8.$$

Summing up, we have shown that $S$ contains an arithmetic progression of primes $p \equiv u\pmod l$.

Then by Theorem \ref{Theorem 1.2},
\begin{equation*}
\left\{ \textrm{prime ideal } \wp \subset \mathcal{O}_{K} : N(\wp)\leq x, N\wp\in S, \pi_{\wp} \mbox{ is onto } \right\}   \gg  \frac{x}{(\log x)^2},
\end{equation*}
Combining Theorem \ref{Theorem 4.1}, we get $[\wp]$ is an Euclidean ideal class.

Finally, we prove that if $h_K=2^t\geq 4$, then $H(K)$ is non-abelian over $\mathbb{Q}$.

Let $\sigma$ be the nontrivial automorphism of $\textrm{Gal}(K/k)$ and $\tilde{\sigma}$ be any lifting of $\sigma$ in $\textrm{Gal}(H(K)/\mathbb{Q})$. The Galois group of Hilbert class field, $\textrm{Gal}(H(K)/K)$ is isomorphic to $\textrm{Cl}(K)$ via Artin reciprocity map. By Chebotarev's density theorem, each element of $\textrm{Gal}(H(K)/K)$ can be represented by Artin symbol $\left(\wp,H(K)/K\right) $ for some prime ideal $\wp$ of $\mathcal{O}_{K}$.

By the property of Artin symbol, we have
\begin{equation}\label{Art1}
\tilde{\sigma}\left(\wp, H(K)/K\right)\tilde{\sigma}^{-1}=\left(\sigma\wp, H(K)/K\right).
\end{equation}
Clearly,
\begin{equation}\label{Art2}
\left(\sigma\wp, H(K)/K\right)\left(\wp, H(K)/K\right)=\left(I\mathcal{O}_{K}, H(K)/K\right)
\end{equation}
where $I$ is an ideal of $k$ such that $I\mathcal{O}_{K}=N_{K/k}(\wp)=\sigma\wp\cdot \wp.$ Since
$I^{h(k)}$ is an principal ideal of $k$, we have
\begin{equation}\label{Art3}
\left(I\mathcal{O}_{K}, H(K)/K\right)^{h(k)}=\textrm{Id}\in \textrm{Gal}(H(K)/K).
\end{equation}
Combining \eqref{Art3} with the assumption $2\nmid h(k)$ and $h(K)=2^t$, we have
\begin{equation}\label{Art4}
\left(I\mathcal{O}_{K}, H(K)/K\right)=\textrm{Id}\in \textrm{Gal}(H(K)/K).
\end{equation}
Combing \eqref{Art1}, \eqref{Art2} and \eqref{Art4}, we get
\begin{equation}\label{Art5}
\tilde{\sigma}\left(\wp, H(K)/K\right)\tilde{\sigma}^{-1}=\left(\wp, H(K)/K\right)^{-1}.
\end{equation}
By assumption $h(K)=2^t\geq 4$. Choose $\wp$ such that $\left(\wp, H(K)/K\right)$ is a generator of $\textrm{Gal}(H(K)/K)$. Then $\left(\wp, H(K)/K\right)^{-1}\neq \left(\wp, H(K)/K\right)$. Combining \eqref{Art5}, we get
\begin{equation*}
\tilde{\sigma}\left(\wp, H(K)/K\right)\tilde{\sigma}^{-1}\neq\left(\wp, H(K)/K\right).
\end{equation*}
Therefore, $\textrm{Gal}(H(K)/\mathbb{Q})$ is a non-abelian group. This concludes the proof.


Based on works of Yue, Bae and  Yue, Ouyang and Zhang, we prove Theorem \ref{genus}, which simplify the condition $ii)$ of Theorem \ref{main1}.

\noindent\textbf{Proof of Theorem \ref{genus}:}

 \emph{Case} 1 follows from Theorem 3.1 of [Bae $\&$ Yue \cite{Yue2011}].

	\emph{Case} 2 follows from Theorem 5.1 of [Bae $\&$ Yue \cite{Yue2011}].

	\emph{Case} 3 follows from Theorem 1.1 and Theorem 1.2 of [Yue \cite{Yue2010}].

	\emph{Case} 4 follows from Theorem 4.1 of [Bae $\&$ Yue \cite{Yue2011}].

	\emph{Case} 5 and \emph{Case} 6 follows from Theorem 3.5 of [Ouyang $\&$ Zhang \cite{Ouyang}].
	
The checkment is rountine. We only illustrate \emph{Case} 5.

In \emph{Case} 5, by Theorem 3.5 of (Ouyang $\&$ Zhang \cite{Ouyang}), $G(K)=\mathbb{Q}(\sqrt{p_{1}},\sqrt{q_{1}},\sqrt{q_{2}})$ if and only if $m=r$, where $q_{1},\cdots,q_{m}$ are all the primes $q_{j} (j=1,2)$ satisfying $\displaystyle\left(\frac{p_{1}}{q_{j}}\right)=1$ and $r$ is the $2$-rank of $\mu_{2}\times\mu_{2}$ generated by $\left(\displaystyle\left(\frac{-1}{q_{j}}\right),\displaystyle\left(\frac{2}{q_{j}}\right)\right) $ with $1\leq j \leq m$ (See equations (3),(4),(5),(6) of [Ouyang $\&$ Zhang \cite{Ouyang}]). We discuss it case by case:\\
If $m=0$, then $r=0$ automatically holds.\\
If $m=1$, then $r=1$ should hold, hence $\left(\displaystyle\left(\frac{-1}{q_{1}}\right),\displaystyle\left(\frac{2}{q_{1}}\right)\right)\neq (1,1) $, which implies $q_{1} \not\equiv 1 \pmod8$.\\
If $m=2$, then $r=2$ should hold, hence $q_{1}, q_{2} \not\equiv 1 \pmod8$ and  $\displaystyle\left(\frac{-1}{q_{1}q_{2}}\right)\neq \displaystyle\left(\frac{2}{q_{1}q_{2}}\right) $. This is equivant to $q_{1} ,q_{2} \not\equiv 1 \pmod8$ and $q_{1}q_{2} \equiv 5,7 \pmod8 $.

Note that in two latter subcases of \emph{Case} 1 and \emph{Case} 2, the original condition of Theorem 3.1 and Theorem 5.1 of \cite{Yue2011} is that $q_1\notin N_{k(\sqrt{\varepsilon})/k}\left(k(\sqrt{\varepsilon})\right)$,  where $N_{k(\sqrt{\varepsilon})/k}$ is the norm map from $k(\sqrt{\varepsilon})$ to $k$. In fact, this is equivalent to the condition
$$
\displaystyle\left(\frac{\varepsilon \mod \mathfrak{q} }{q_{1}}\right)=-1$$
where $\mathfrak{q}$ is some (equivalently any) prime ideal of $\mathcal{O}_k$ lying above $q_1$.

Indeed, by Hasse norm theorem and properties of Hilbert symbol,
$$q_1\notin N_{k(\sqrt{\varepsilon})/k}\left(k(\sqrt{\varepsilon})\right)\Leftrightarrow
\displaystyle\left(\frac{\varepsilon, q_{1} }{\mathfrak{p}}\right)=-1\ \mathrm{for\ some\ prime\ ideal}\  \mathfrak{p}|2q_1\mathcal{O}_k.
$$

Since $q_1\equiv 1\pmod{4}$ splits in $k=\mathbb{Q}(\sqrt{p_1})$, we have
\begin{equation}\label{f1}\displaystyle\left(\frac{\varepsilon, q_{1} }{\mathfrak{q}}\right)=\displaystyle\left(\frac{\varepsilon \mod \mathfrak{q} }{q_{1}}\right),\ \mathrm{where}\ \mathfrak{q}|q_1\mathcal{O}_k.\end{equation}
This implies that
\begin{equation}\label{f4}\prod_{\mathfrak{q}|q_1\mathcal{O}_k}\displaystyle\left(\frac{\varepsilon, q_{1} }{\mathfrak{q}}\right)=\left(\frac{N(\varepsilon)}{q_1}\right)=\left(\frac{-1}{q_1}\right)=1
\end{equation}
since the norm of the fundamental unit, $N(\varepsilon)=-1$ for general quadratic field $k=\mathbb{Q}(\sqrt{p_1})$ with prime $p_1\equiv 1\pmod{4}$.

If $p_1\equiv 5\pmod{8}$, then $2$ inerts in $k=\mathbb{Q}(\sqrt{p_1})$. By the product formula of Hilbert symbol and \eqref{f4}, we get
\begin{equation}\label{f2}\displaystyle\left(\frac{\varepsilon, q_{1} }{2\mathcal{O}_k}\right)=1.\end{equation}
Otherwise, if $p_1\equiv 1\pmod{8}$, then $2$ splits in $k=\mathbb{Q}(\sqrt{p_1})$. Then we have
\begin{equation}\label{f3}\displaystyle\left(\frac{\varepsilon, q_{1} }{\mathfrak{p}}\right)=(-1)^{\frac{(\varepsilon-1)(q_1-1)}{4}}=1,\ \mathrm{where}\ \mathfrak{p}|2\mathcal{O}_k,\end{equation}
where $\varepsilon$ is viewed as an element of the $2$-adic field $\mathbb{Q}_2$.

Combining \eqref{f1}, \eqref{f4}, \eqref{f2} and \eqref{f3}, we get the desired result.

\noindent\textbf{Proof of Theorem \ref{main3}:}
If $G(K)=\mathbb{Q}(\sqrt{p_{1}},\sqrt{q_{1}},\sqrt{q_{2}})$, then by Theorem \ref{main1}, $K$ has an Euclidean ideal.

Now assume that $K$ has an Euclidean ideal. Then $\mathrm{Cl}(K)$ is a cyclic group. By assumption $h(K)=2^t$. Therefore
\begin{equation}\label{Converse1}
\textrm{Cl}(K)/\textrm{Cl}(K)^{2}\cong \mathbb{Z}/2\mathbb{Z}.
\end{equation}
which implies that $G(K)$ is a quadrati extension of $K$.

Since $\mathbb{Q}(\sqrt{p_{1}},\sqrt{q_{1}},\sqrt{q_{2}})=K\mathbb{Q}(\sqrt{q_1})$ and $q_1\equiv 1 \pmod{4}$, we know that $\mathbb{Q}(\sqrt{p_{1}},\sqrt{q_{1}},\sqrt{q_{2}})/K$ is unramified outside $q_1$. Similarly, since $\mathbb{Q}(\sqrt{p_{1}},\sqrt{q_{1}},\sqrt{q_{2}})=K\mathbb{Q}(\sqrt{q_2})$ and $q_2\equiv 1 \pmod{4}$, $\mathbb{Q}(\sqrt{p_{1}},\sqrt{q_{1}},\sqrt{q_{2}})/K$ is also unramified outside $q_2$. Therefore $\mathbb{Q}(\sqrt{p_{1}},\sqrt{q_{1}},\sqrt{q_{2}})/K$ is an unramified extension. By definition of Hilbert genus field, $\mathbb{Q}(\sqrt{p_{1}},\sqrt{q_{1}},\sqrt{q_{2}})\subset G(K)$. Comparing the degrees, we get
$G(K)=\mathbb{Q}(\sqrt{p_{1}},\sqrt{q_{1}},\sqrt{q_{2}})$.

The rest part follows by checking the conditions of \emph{Cases} 1, 4, 5 in Theorem \ref{genus} under the assumption $q_1\equiv q_2\equiv 1\pmod{4}$.
\section{Open Questions}
In this section, we will list some open questions concerning the existence of Euclidean ideals.

\textbf{Q1.} Remove the condition: $q_1\neq 3$ and $q_2\neq3$ in Theorem \ref{main1}.

Under the conditions $G(K)=\mathbb{Q}(\sqrt{p_{1}},\sqrt{q_{1}},\sqrt{q_{2}})$ and $h_K=2^t$ in Theorem \ref{main1}, $\textrm{Cl}(K)$ is a cyclic group of order $2^t$. Assuming GRH, $K$ should have an Euclidean ideal by Lenstra's result. The unnatural restriction $q_1\neq 3, \ q_2\neq3$ comes from the use of Theorem \ref{Theorem 1.2}, which relies on the deep analytic result of \cite{Heath, Narkiewicz}. At the present, we have no idea on how to improve it.

\textbf{Q2.} Using the theory of Hilbert genus fields, find all the biquadratic fields $K$ of the form $K=\mathbb{Q}(\sqrt{p_1}, \sqrt{q_1q_2})$ having Euclidean ideals, under the condition $h_K=2^t\ (t\geq1)$, where $p_1, q_1, q_2$ are distinct prime numbers.

Under the condition $h_K=2^t\ (t\geq1)$, $\textrm{Cl}(K)$ is cyclic if and only if $G(K)$ is a quadratic extension of $K$. For $K=\mathbb{Q}(\sqrt{p_1}, \sqrt{q_1q_2})$, $G(K)$ is not always equal to $\mathbb{Q}(\sqrt{p_{1}},\sqrt{q_{1}},\sqrt{q_{2}})$ from Theorem \ref{genus} and the works of Yue \cite{Yue2010}, Bae $\&$ Yue \cite{Yue2011}, Ouyang $\&$ Zhang \cite{Ouyang}. So \textbf{Q2} is amount to relax the condition $G(K)=\mathbb{Q}(\sqrt{p_{1}},\sqrt{q_{1}},\sqrt{q_{2}})$ to $[G(K):K]=2$ in Theorem \ref{main1}.

\textbf{Q3.} Using the theory of Hilbert genus fields, find all the biquadratic fields $K$ of the form $K=\mathbb{Q}(\sqrt{p_1}, \sqrt{d})$ having Euclidean ideals, under the condition $h_K=2^t\ (t\geq1)$, where $p_1$ is a prime and $d>0$ is a square-free integer.

This is a generalization of \textbf{Q2}. Under the condition $h_K=2^t\ (t\geq1)$, $\textrm{Cl}(K)$ is cyclic if and only if $r_2(\textrm{Cl}(K))=1$. By formula \eqref{2dim}, $\omega(d)\leq 5$, where $\omega(d)$ is the number of distinct prime factors of $d$.

\textbf{Q4.} Using the theory of Hilbert genus fields, find all the biquadratic fields $K$ of the form $K=\mathbb{Q}(\sqrt{\delta}, \sqrt{d})$ having Euclidean ideals, under the condition $h_K=2^t\ (t\geq1)$, where $\mathbb{Q}(\sqrt{\delta})$ is a real quadratic field with odd class number, and $d>0$ is square-free.

This is a further generalization of \textbf{Q3}. Note that $\mathbb{Q}(\sqrt{\delta})$ has odd class number if and only if $\delta=p, 2p_1$ or $p_{1}p_{2}$
where $p, p_{1}$ and $p_{2}$ are primes with $p_1\equiv p_2\equiv 3\pmod{4}$. The works of Yue, Bae and Yue, Ouyang and Zhang
completed the construction of the Hilbert genus field of all real biquadartic fields $K=\mathbb{Q}(\sqrt{\delta}, \sqrt{d})$.

\textbf{Q5.} Find an infinitely family of biquadratic fields $K$ having Euclidean ideals.

If one looks for such family in Theorem \ref{main1}, one may find that the difficulty is to find infinitely many prime triples $(p_1, q_1, q_2)$  whose class numbers satisfy
$$h_{\mathbb{Q}(\sqrt{p_{1}})}=1;\  h_{\mathbb{Q}(\sqrt{q_{1}q_{2}})}\ and \ h_{\mathbb{Q}(\sqrt{p_{1}{q_{1}{q_{2}}}})}\ are\  powers\ of\ 2.$$

\textbf{Q6.} Find an example of real biquadratic field $K$ having Euclidean ideals without genus theory.

  This is amount to find an example of real biquadratic field $K$ having Euclidean ideals, whose class number is not a power of $2$.

\textbf{Acknowledgement:} We are grateful to Nick Rome for pointing out a mistake of earlier version and sending some enlightening suggestion to us, which made the paper more elegant and readable.
This work is supported by the Laboratory Open Project funded by the State Key Laboratory of Cryptology in Beijing.

\bibliography{central}

\end{document}